\documentclass[a4paper,12pt,reqno]{amsart}
\usepackage[headings]{fullpage}

\usepackage{hyperref}

\usepackage{amsthm}
\usepackage{amsmath}
\usepackage{amssymb}

\usepackage[all]{xy}   

\theoremstyle{plain}
\newtheorem{thm}{Theorem}[section]
\newtheorem*{thm*}{Theorem}

\newtheorem{lem}[thm]{Lemma}
\newtheorem{cor}[thm]{Corollary}

\theoremstyle{definition}
\newtheorem{defn}[thm]{Definition}

\theoremstyle{remark}
\newtheorem*{rem*}{Remark}
\newtheorem*{notation*}{Notation}
\newtheorem*{ex*}{Example}
\newtheorem*{ackn*}{Acknowledgements}

\newcommand{\Vol}{\operatorname{Vol}}
\newcommand{\Id}{\operatorname{Id}}

\def\co{\colon\thinspace}
\begin{document}

\title[Einstein metrics, volume entropy, and exotic smooth structures]{An essential relation between Einstein metrics, volume entropy, and exotic smooth structures}

\author{Michael Brunnbauer}
\address{Mathematisches Institut, Ludwig-Maximilians-Universit\"at M\"unchen, Theresienstr.\ 39, D-80333 M\"unchen, Germany}
\email{michael.brunnbauer@mathematik.uni-muenchen.de}

\author{Masashi Ishida}
\address{Department of Mathematics, Sophia University, 7-1 Kioi Cho, Chiyoda-Ku, Tokyo 102-8554, Japan}
\email{masashi@math.sunysb.edu}

\author{Pablo Su\'arez-Serrato}
\address{Mathematisches Institut, Ludwig-Maximilians-Universit\"at M\"unchen, Theresienstr.\ 39, D-80333 M\"unchen, Germany}
\email{p.suarez-serrato@cantab.net}

\date{\today}


\begin{abstract}
We show that the minimal volume entropy of closed manifolds remains unaffected when nonessential manifolds are added in a connected sum. We combine this result with the stable cohomotopy invariant of Bauer--Furuta in order to present an infinite family of four--manifolds with the following properties:

 \begin{enumerate}
\item They have positive minimal volume entropy.
\item They satisfy a strict version of the Gromov--Hitchin--Thorpe inequality, with a minimal volume entropy term.
 \item They nevertheless admit infinitely many distinct smooth structures for which no compatible Einstein metric exists.
 \end{enumerate}
\end{abstract}

\maketitle


\section{Introduction}

The volume of a ball of radius $R$ in the universal Riemannian covering $(\widetilde{M},\tilde{g})$ of a Riemannian manifold $(M,g)$ grows at an exponential rate
\[ \lambda (M,g):=\lim\limits_{R \to \infty} R^{-1}\log \, \Vol (B_{\tilde g}(x,R)) \geq 0. \]
This geometric invariant is called the \emph{volume entropy} of $(M,g)$. Manning \cite{Manning(1979)} showed that this limit exists and does not depend on the choice of the center point $x\in\widetilde{M}$.

The \emph{minimal volume entropy} $\lambda(M)$ is by definition the infimum of $\lambda (M,g)$ as it runs over all smooth unit volume metrics of $M$. The value $\lambda(M)$ of this topological invariant has been shown by the first author to depend only on the image of the fundamental class of $M$ under the classifying map of the universal covering \cite{Brunnbauer(2007a)}. 

A \emph{word metric} on a group $G$ is a way to measure distance between two elements of $G$. A choice of a set $S$ which generates $G$ induces a distance function $d_{S}(g,h)$ between elements $g$ and $h$ in $G$, defined as the minimum number of elements of $S\cup S^{-1}$ needed to express $g^{-1}h$. The growth function of a group $G$ with respect to a symmetric generating set $S$ describes the size of balls in $G$, that is, it counts the number of elements of $G$ that can be written as a product of $n$ elements of $S$.

Let $B_{n}(G,S):=\{ x \in G \, | \,  d_{S}(x,e) \leq n \}$ be the ball of radius $n$. The growth function of $G$ is defined as $\beta (n):=|B_{n}(G,S)|$. If $\beta (n) \leq C(n^{k}+1)$ for some constant $C$ and $k < \infty$, we say $G$ has polynomial growth rate. If  $\beta (n) \geq a^{k}$ for some $a>1$, we say $G$ has exponential growth rate. The type of growth rate of the group $G$ is independent of the choice of $S$.

A notorious property states that if the growth of the fundamental group $\pi_1(M)$ is of exponential type then $\lambda (M,g) > 0$ for every smooth metric $g$ on $M$. A subtle point to notice is that the infimum $\lambda(M)$ may vanish even in spite of the growth of $\pi_1(M)$ being exponential.

Consider a pair of closed manifolds $M$ and $N$. Even if we assume we know the values of $\lambda(M)$ and $\lambda(N)$, the question of how to determine $\lambda(M \# N)$ remains open in general. This is due to the fact that the fundamental group (and thus the universal covering) of $M\#N$ may be `large' in comparison to the fundamental groups of $M$ and $N$. Indeed, if $A$ and $B$ are two finitely generated groups, then the free product $A \ast B$ contains a free subgroup of rank two---and thus grows exponentially---unless $A$ is trivial or $B$ is trivial, or $A$ and $B$ are both of order two. The first main result of the present paper is to compute the minimal volume entropy of $M\# N$ with some topological assumptions on the second summand.

Let $\Psi\co N\to K(\pi_1(N),1)$ denote the classifying map of the universal covering. The homotopy class of this map is uniquely defined, and if one may choose $\Psi$ to have image in the $(n-1)$--skeleton of the Eilenberg--Mac\,Lane space, then $N$ is called \emph{nonessential}. This property allows us to show:

\begin{thm}\label{main-A}
Let $M$ and $N$ be two connected closed manifolds. If $N$ is orientable and nonessential, then
\[ \lambda(M\# N) = \lambda(M). \]
\end{thm}

We will use this result to extract new information about exotic smooth structures on four--manifolds and existence of Einstein metrics. Recall that any Riemannian metric is called Einstein if its Ricci curvature, considered as a function on the unit tangent bundle, is constant. There exist closed manifolds that do not admit any Einstein metric at all. It is known that any closed Einstein four--manifold $X$ must satisfy the following topological constraint
\begin{eqnarray*}\label{ht-int}
2\chi(X) \geq 3|\tau(X)|,
\end{eqnarray*}
where $\chi(X)$ and $\tau(X)$ denote respectively the Euler characteristic and the signature of $M$. This inequality is called the Hitchin--Thorpe inequality \cite{hit, thor}. In particular, Hitchin \cite{hit} proved that any closed oriented Einstein four--manifold satisfying $2\chi(X) = 3|\tau(X)|$ is finitely covered by either $K3$ surface or 4--torus.

On the other hand, using the Seiberg--Witten monopole equations \cite{w}, LeBrun \cite{leb-1996} found a new obstruction to the existence of Einstein metrics in dimension four. This obstruction provided the first means of exhibiting the dependence of the existence of Einstein metrics on smooth structures of underlying four--manifolds. By using the new obstruction, examples of pairs of homeomorphic, but not diffeomorphic, simply-connected four--manifolds such that one manifold admits an Einstein metric and satisfies the strict Hitchin--Thorpe inequality while the other cannot admit any Einstein metric, were first found by Kotschick \cite{kot(1998)}. Further examples can also be consulted in \cite{leb-2001}.

The Hitchin--Thorpe inequality can be strengthened with a simplicial volume term to
\begin{eqnarray*}\label{ght-int}
2\chi(X) \geq 3|\tau(X)|+\frac{1}{81\pi^2}\|X\|,
\end{eqnarray*}
which is called the Gromov--Hitchin--Thorpe inequality \cite{kot(1998-1)}.

Moreover, Kotschick \cite{kot(2004)} improved the Gromov--Hitchin--Thorpe inequality by introducing a minimal volume entropy term. In fact, he pointed out that any closed Einstein four--manifold $X$ must satisfy the following inequality
\begin{eqnarray*}\label{gro-vol}
2\chi(X) - 3|\tau(X)| \geq \frac{1}{54{\pi}^2}{\lambda}(X)^4,
\end{eqnarray*}
where equality can occur if and only if every Einstein metric on $X$ is flat, is a non-flat Calabi--Yau metric, or is a metric of constant negative sectional curvature. In this article, we shall call the strict case of the previous inequality the strict Gromov--Hitchin--Thorpe inequality, that is,
\begin{eqnarray}\label{ght}
2\chi(X) - 3|\tau(X)| > \frac{1}{54{\pi}^2}{\lambda}(X)^4.
\end{eqnarray}

Although related ideas and results appear in the work of Besson--Courtois--Gallot \cite[p. 774]{BCG(1995)} and Sambusetti \cite[p. 537]{s}, the above inequality is found in the literature only in \cite{kot(2004)}.

This strengthened version of the Gromov--Hitchin--Thorpe inequality contains all known homotopy invariant obstructions---previously found by Hitchin \cite{hit}, Thorpe \cite{thor}, Gromov \cite{gromov(1982)}, Besson--Courtois--Gallot \cite{BCG(1995)}, Sambusetti \cite{s} and Kotschick\cite{kot(1998-1)}---to the existence of Einstein metrics in dimension four as special cases \cite{kot(2004)}.

To the best of our knowledge it is still an open question if there exists a topological 4-manifold which has positive minimal volume entropy, satisfies the strict Gromov--Hitchin--Thorpe inequality (\ref{ght}), and nevertheless admits infinitely many distinct smooth structures for which no compatible Einstein metric exists (cf.\ \cite{IS(2008)}). In this article, however, we are able to give an affirmative answer to this question as follows.

\begin{thm}\label{main-B}
There exists an infinite family of topological four--manifolds which have positive minimal volume entropy, satisfy the strict Gromov--Hitchin--Thorpe inequality (\ref{ght}) and nevertheless admit infinitely many distinct smooth structures for which no compatible Einstein metric exists.
\end{thm}

In section \ref{min_entr_conn_sum} the behaviour of the minimal volume entropy under connected sum is considered in some special cases. This result is then used in the last section to prove Theorem \ref{main-B}, where we briefly review some of the Seiberg--Witten theory that is employed.

\section*{Acknowledgements}

 Michael Brunnbauer and Pablo Su\'arez-Serrato are grateful to the DFG (Deutsche  Forschungsgemeinschaft) for financial support. Masashi Ishida is partially supported by the Grant-in-Aid for Scientific Research (C), Japan Society for the Promotion of Science, No. 20540090.


\section{Minimal volume entropy and connected sums}\label{min_entr_conn_sum}

There is not much known about the behaviour of the minimal volume entropy under connected sum. In fact, it seems that the only result (in arbitrary dimensions) was shown by Babenko in \cite{Babenko(1995)}, Corollary 2. He proved that adding a simply-connected summand does not change the minimal volume entropy (see also \cite{Brunnbauer(2007a)}, Corollary 8.3). We will extend this result to orientable manifolds which are not essential in the following sense:

\begin{defn}
A connected closed manifold $M$ of dimension $n$ is called \emph{essential} if there exists a map $M\to K$ to an aspherical complex $K$ that does not contract to the $(n-1)$--skeleton of $K$.
\end{defn}

This notion was introduced by Gromov in \cite{gromov(1983)}. Recall that the homotopy class of maps to aspherical complexes is uniquely determined by the induced homomorphism on fundamental groups. Therefore, it suffices to consider the classifying map of the universal covering $M\to K(\pi_1(M),1)$ since every map $M\to K$ factorizes over this map. Thus, for manifolds the definition of essentialness implies that the classifying map contracts to the $(n-1)$--skeleton if and only if $M$ is nonessential.

In particular, every simply-connected manifold is nonessential since the respective classifying space is a point. Further examples of nonessential manifolds are provided by products of arbitrary manifolds with simply-connected manifolds (of positive dimension) since the classifying map factors over the projection map to the first factor of the product.

Note that nonessential manifolds have zero minimal volume entropy (see \cite{Babenko(1992)}, Proposition 2.7). We will prove that adding orientable nonessential summands does not change the minimal volume entropy (Theorem \ref{main-A}).

\begin{rem*}
Note that the fundamental group of $M$ may change dramatically if one adds a summand $N$ that is not simply-connected. For instance, if $\pi_1(M)$ is of polynomial growth, as we mentioned in the introduction it may happen that $\pi_1(M\#N) = \pi_1(M)*\pi_1(N)$ is of exponential growth. Thus the universal covering of $M\# N$ is enormous compared to the universal covering of $M$.
\end{rem*}

Theorem \ref{main-A} allows us to improve \cite{IS(2008)}, Lemma 69 (for $k=1$ in their notation).

\begin{cor}\label{main-cor}
Let $X_m,m=1,\ldots,n,$ be closed simply-connected $4$-manifolds. Consider the connected sum
\[ M := (\#_{m=1}^n X_m) \# (\Sigma_g \times \Sigma_h) \# \ell_1(S^1 \times S^3) \# \ell_2\overline{\mathbb{C}\mathrm{P}^2}, \]
where $g,h \geq 1$ and $\ell_1,\ell_2\geq 0$. Then the minimal volume entropy of $M$ satisfies the following bounds
\[ 16 (g-1)(h-1) \leq \lambda(M)^4 \leq 256\pi^2 (g-1)(h-1). \]
\end{cor}

\begin{proof}
Note that the $X_m$ and $\overline{\mathbb{C}\mathrm{P}^2}$ are simply-connected and that $S^1\times S^3$ is nonessential. Thus, $\lambda(M) = \lambda(\Sigma_g\times\Sigma_h)$ by Theorem \ref{main-A}. By \cite{Babenko(1992)}, Proposition 2.6
\begin{align*}
\lambda(\Sigma_g\times\Sigma_h)^4 &\leq 256 \cdot \lambda(\Sigma_g)^2/4 \cdot \lambda(\Sigma_h)^2/4 \\
&= 256\pi^2 (g-1)(h-1)
\end{align*}
since the minimal volume entropy fulfills $\lambda(\Sigma_g)^2 = 4\pi(g-1)$ for $g\geq 1$.

The lower bound stems from the inequality
\[ {\textstyle\frac{n^{n/2}}{n!}} \|N\| \leq \lambda(N)^n \]
for arbitrary connected closed orientable manifolds $N$ of dimension $n$ (see Besson--Courtois--Gallot \cite{BCG(1991)}, Th\'eo\-r\`emes 3.8 and 3.16) and from Bucher-Karlsson's computation of the simplicial volume
\[ \|\Sigma_g\times \Sigma_h\| = 24(g-1)(h-1) \]
(see \cite{Bucher-Karlsson(2007)}, Corollary 3).
\end{proof}

We will split the proof of Theorem \ref{main-A} into two lemmas. But first note that it suffices to consider the case $n\geq 3$ since the only surface that is not essential is the sphere.

\begin{lem}
Let $M$ and $N$ be two connected closed manifolds of dimension $n\geq 3$. If $N$ is orientable, then $\lambda(M\# N) \geq \lambda(M)$.
\end{lem}

\begin{proof}
The map $f \co M\#N \to M$ that contracts $N$ to one point is orientation-true (i.\,e.\ it maps orientation preserving loops to orientation preserving ones and orientation reversing loops to orientation reversing ones) and of absolute degree one. For more details on the notion of absolute degree see for example section 2 of \cite{Brunnbauer(2007a)}. However, note that for maps between orientable manifolds the absolute degree is just the absolute value of the usual mapping degree.

It is shown in \cite{Brunnbauer(2007a)}, Corollary 8.2 that in this situation the minimal volume entropy of $M$ equals the `relative' minimal volume entropy of $M\# N$ with respect to the epimorphism $f_* \co \pi_1(M\#N)=\pi_1(M)*\pi_1(N) \twoheadrightarrow \pi_1(M)$. This `relative' invariant is defined as follows:
\[ \lambda_{f_*} (M\# N) := \inf_g \lambda_{f_*} (M\# N, g) \Vol(M\#N,g)^{1/n} \]
where the infimum is taken over all Riemannian metrics $g$ on $M\# N$ and $\lambda_{f_*} (M\# N, g)$ is the volume entropy of the covering $\overline{(M\#N)}_{f_*}$ of $M\# N$ associated to the subgroup $\ker f_*\subset\pi_1(M\# N)$ with respect to the lifted metric. Note that this covering may be obtained from the universal covering of $M$ by pullback along $f$.

So, we have $\lambda_{f_*}(M\# N) = \lambda(M)$ by \cite{Brunnbauer(2007a)}, Corollary 8.2. Since a ball in the universal covering has bigger volume than a ball of the same radius in the covering $\overline{(M\#N)}_{f_*}$, it follows from the definition of volume entropy that $\lambda(M\# N)\geq \lambda_{f_*}(M\# N)$. Therefore, $\lambda(M\# N)\geq \lambda(M)$.
\end{proof}

To show the converse inequality, we have to extend the notion of minimal volume entropy to connected finite simplicial complexes by using continuous piecewise smooth Riemannian metrics (see \cite{Babenko(1992), Sabourau(2006)} and also \cite{Brunnbauer(2007a)} for more details).

\begin{defn}
Consider a connected finite simplicial complex $X$ of dimension $n$. An \emph{extension} of $X$ is a simplicial complex $X'$ that is obtained from $X$ by successive attachment of finitely many cells of dimension strictly less than $n$ such that the inclusion $X\subset X'$ induces an injective homomorphism on fundamental groups $\pi_1(X)\hookrightarrow \pi_1(X')$.
\end{defn}

In \cite{Brunnbauer(2007a)}, section 9 (in particular Theorem 9.5) the following theorem was derived:

\begin{thm}\label{ext_ax}
If $X'$ is an extension of $X$, then $\lambda(X')=\lambda(X)$.
\end{thm}

Using this result, it is not hard to show the remaining inequality and to complete the proof of Theorem \ref{main-A}.

\begin{lem}
Let $M$ and $N$ be two connected closed manifolds of dimension $n\geq 3$. If $N$ is nonessential, then $\lambda(M\# N) \leq \lambda(M)$.
\end{lem}

\begin{proof}
Let $\Psi\co N\to K(\pi_1(N),1)$ be the classifying map of the universal covering. Since $N$ is nonessential we may assume that the image of $\Psi$ lies in the $(n-1)$--skeleton and by compactness that it lies in a finite subcomplex $K\subset K(\pi_1(N),1)^{(n-1)}$. Denote this map by $f \co N\to K$.

The composition
\[ M\# N \to M\vee N \xrightarrow{\Id\vee f} M\vee K, \]
where the first map contracts the belt sphere of the connected sum to a point, induces an isomorphism on fundamental groups and is $(n,1)$--monotone, i.\,e.\ the preimage of an open $n$--simplex consists of at most one open $n$--simplex. Hence, $\lambda(M\# N) \leq \lambda(M\vee K)$ by \cite{Sabourau(2006)}, Lemma 3.5 (see also \cite{Brunnbauer(2007a)}, Lemma 4.1).

But $M\vee K$ is an extension of $M$. Therefore, $\lambda(M\vee K) = \lambda(M)$ by Theorem \ref{ext_ax}. This proves the lemma and thus Theorem \ref{main-A}, too.
\end{proof}


\section{Einstein metrics and exotic smooth structures}

As we have already mentioned in the introduction, LeBrun \cite{leb-1996} found a new obstruction to the existence of Einstein metrics on four--manifolds. Armed with this new obstruction, Kotschick \cite{kot(1998)} found the first examples of pairs of homeomorphic, but not diffeomorphic, simply-connected four--manifolds such that one manifold admits an Einstein metric and satisfies the strict Hitchin--Thorpe inequality while the other cannot admit any Einstein metric.

The technique relies on the existence of solutions to the Seiberg--Witten monopole equations \cite{w}, i.e., the existence of monopole classes \cite{kro}. Let $X$ be a closed smooth four--manifold $X$ with $b^{+}(X)>1$, where ${b}^+(X)$ denotes the dimension of the maximal positive definite linear subspace in the second cohomology of $X$. The first Chern class $c_{1}({\mathcal L}_{\mathfrak{s}})$ of complex line bundle ${\mathcal L}_{\mathfrak{s}}$ associated with a spin${}^c$--structure $\mathfrak{s}$ of $X$ with $b^{+}(X)>1$ is called monopole class of $X$ if the corresponding Seiberg--Witten monopole equations have a solution for every choice of Riemannian metric \cite{kro}. It is known that the set of all monopole classes of  $X$ is a finite set \cite{ishi-leb-2}.

There are some methods to detect the existence of monopole classes. The first one is due to Witten, who introduced in the celebrated article \cite{w} an invariant of smooth 4-manifolds using the fundamental homology class of the Seiberg--Witten moduli space. This is called the Seiberg--Witten invariant, and is well-defined for any closed four--manifold $X$ with $b^{+}(X)>1$. The Seiberg--Witten invariant defines an integer valued function $SW_{X}$ over the set of all isomorphism classes of spin${}^c$ structures of $X$ with $b^{+}(X)>1$. The non-triviality of $SW_{X}$ for some spin${}^c$--structure $\mathfrak{s}$ tells us that $c_{1}({\mathcal L}_{\mathfrak{s}})$ is a monopole class of $X$.

The second one is due to Bauer and Furuta \cite{b-f, b-1}. It is a sophisticated refinement of the idea of the construction of $SW_{X}$. It detects the presence of monopole classes by an element of a certain complicated stable cohomotopy group ${\pi}^{{b}^+}_{S^1, \mathcal{B}}({\rm Pic}^0(X), {\rm ind} D)$, where ${b}^+:={b}^+(X)$  and ${\rm ind} D$ is the virtual index bundle for the Dirac operators parametrised by the    $b_{1}(X)$--dimensional Picard torus ${\rm Pic}^0(X)$. This invariant is called the stable cohomotopy Seiberg--Wittten invariant, denote it by $BF_{X}$:
\begin{eqnarray*}\label{b-f-inv}
BF_{X}(\mathfrak{s}) \in {\pi}^{{b}^+}_{S^1, \mathcal{B}}({\rm Pic}^0(X), {\rm ind} D).
\end{eqnarray*}
Bauer \cite{b-1} proved a non-vanishing theorem for $BF_{X}$ which showed that $BF_{X}$ is strictly stronger than $SW_{X}$. It is also known that the non-triviality of $BF_{X}$ for a spin${}^c$--structure $\mathfrak{s}$ ensures us that $c_{1}({\mathcal L}_{\mathfrak{s}})$ is a monopole class of $X$ (cf. \cite{ishi-leb-2}). \par
By showing a new non-vanishing theorem for $BF_{X}$ and combining it with the technique of LeBrun, the following obstruction to the existence of Einstein metrics in dimension four was derived in \cite{IS(2008)}, Corollary 70.
\begin{thm}\label{speical-ein}
For $m=1,2,3$, let $X_{m}$ be simply connected symplectic four--manifolds with ${b}^{+}(X_{m}) \equiv 3 \ (\bmod \ 4)$. Consider the following  connected sum
\begin{eqnarray*}
M:=(\#_{m=1}^{n} X_{m}) \# k (\Sigma_{h} \times \Sigma_{g}) \# \ell_{1}({S}^{1} \times {S}^{3}) \# \ell_{2} \overline{{\mathbb C}{P}^{2}},
\end{eqnarray*}
where $n, k \geq 1$ satisfying $n+k \leq 3$, $\ell_{1}, \ell_{2} \geq 0$ and $g, h$ are odd integers $\geq 1$. Then $M$ cannot admit any Einstein metric if
\begin{eqnarray*}
4(n+\ell_{1} + k) + \ell_{2} \geq \frac{1}{3}\Big( \sum_{m=1}^{n} 2\chi(X_{m})+3\tau(X_{m})+4k(1-h)(1-g) \Big).
\end{eqnarray*}
\end{thm}
We shall use this theorem to prove Theorem \ref{main-B}. \par
On the other hand, let us recall a construction of a certain sequence of homotopy $K3$ surfaces. See also \cite{BPV}. Let $Y_{0}$ be a Kummer surface with an elliptic fibration $Y_{0} \rightarrow {\mathbb C}{P}^{1}$. Let $Y_{\ell}$ be obtained from $Y_{0}$ by performing a logarithmic transformation of order $2 \ell + 1$ on a non-singular fiber of $Y_{0}$. Then, $Y_{\ell}$ are simply connected spin manifolds with $b^{+}(Y_{\ell}) = 3$ and $b^{-}(Y_{\ell}) = 19$. By the Freedman classification \cite{freedman}, $Y_{\ell}$ must be homeomorphic to a $K3$ surface. And $Y_{\ell}$ is a K{\"{a}}hler surface with $b^{+}(Y_{\ell}) > 1$ and hence a result of Witten \cite{w} tells us that $\pm {c}_{1}(Y_{\ell})$ are monopole classes of $Y_{\ell}$ for each $\ell$.

Moreover, we need to recall a result of Gompf \cite{gom}, who showed that for arbitrary integers $\alpha \geq 2$ and $\beta \geq 0$, there exists a simply connected symplectic spin 4--manifold $X_{\alpha, \beta}$ satisfying $\chi(X_{\alpha, \beta}) = 24\alpha+4\beta$ and $\tau(X_{\alpha, \beta}) = -16\alpha$. So we obtain:
\begin{eqnarray*}
 b^+(X_{\alpha, \beta}) &=& 4\alpha+2\beta-1 \\
 2\chi(X_{\alpha, \beta}) + 3\tau(X_{\alpha, \beta}) &=& 8\beta \\
 2\chi(X_{\alpha, \beta}) - 3\tau(X_{\alpha, \beta}) &=& 8(12\alpha+\beta)
\end{eqnarray*}

A celebrated result of Taubes \cite{t-1} tells us that the manifolds $X_{\alpha, \beta}$ that Gompf considered have $\pm {c}_{1}(X_{\alpha, \beta})$ as monopole classes.

With these constructions, Corollary \ref{main-cor} and Theorem \ref{speical-ein} enable us to prove the following result.
\begin{thm}\label{volume-entropy-1}
There exists an infinite family of topological spin four--manifolds satisfying the following three properties:
\begin{itemize}
\item Each 4--manifold $X$ has positive volume entropy, ${\lambda}(X)>0$.
\item Each 4--manifold $X$ satisfies the strict Gromov--Hitchin--Thorpe inequality (\ref{ght}),
\begin{eqnarray*}
2\chi(X) - 3|\tau(X)| > \frac{1}{54{\pi}^2}{\lambda}(X)^4.
\end{eqnarray*}
 \item Each 4--manifold $X$ admits infinitely many distinct smooth structures for which no compatible Einstein metric exists.
\end{itemize}
\end{thm}

\begin{proof}
For any fixed pair $(g, h)$ of odd integers $\geq 3$, by taking a suitable pair $(m, n)$ of integers satisfying $4m+2n-1 \equiv 3 \ (\bmod \ 4)$, $m \geq 2$ and $n \geq 1$, we are always able to find at east one positive integer $\ell_{1}$ satisfying the following three inequalities
\begin{eqnarray}\label{ein-in-1}
2n-\frac{5}{27}(g-1)(h-1)-3 > \ell_{1}.
\end{eqnarray}
\begin{eqnarray}\label{ein-in-2}
2(n+12m)-\frac{5}{27} (g-1)(h-1)+21 > \ell_{1}.
\end{eqnarray}
\begin{eqnarray}\label{ein-in-3}
\ell_{1} \geq \frac{1}{3}\Big( 2n+(g-1)(h-1) \Big)-3.
\end{eqnarray}
Notice that we have infinitely many choices of such pairs $(m, n)$ and, hence, of $\ell_{1}$.  For each such five integers $(m, n, g, h, \ell_{1})$ and for any new positive integer $\ell$, consider the following connected sum of four--manifolds
\begin{eqnarray*}
M(m, n, \ell, g, h, \ell_{1}):=X_{m, n} \# {Y}_{\ell} \# (\Sigma_{g} \times \Sigma_{h}) \# \ell_{1}({S}^{1} \times {S}^{3}),
\end{eqnarray*}
where $X_{m, n}$ is Gompf's simply connected symplectic spin 4-manifold \cite{gom} and ${Y}_{\ell}$ is obtained from $Y_{0}$ by performing a logarithmic transformation of order $2 \ell + 1$ on a non-singular fiber of $Y_{0}$.  \par
It is clear that $M(m, n, \ell, g, h, \ell_{1})$ is homeomorphic to the following spin 4--manifold
\begin{eqnarray}\label{homeo-spin}
X_{m, n} \# K3 \# (\Sigma_{g} \times \Sigma_{h}) \# \ell_{1}({S}^{1} \times {S}^{3}).
\end{eqnarray}
For any fixed $(m, n, g, h, \ell_{1})$, the sequence $\{M(m, n, \ell, g, h, \ell_{1}) \ | \ \ell \in {\mathbb N} \}$ contains infinitely many distinct diffeotypes. One can use the bandwidth argument (cf. \cite{ishi-leb-1, ishi-leb-2}) to see this. Alternatively, one can also see this by using only the finiteness property of the set of special monopole classes (cf. \cite{kot(2004)} and compare with the results in \cite{kot(2004-2)}). \par
Moreover, we will show that $M(m, n, \ell, g, h, \ell_{1})$ cannot admit any Einstein metric. First of all, notice that $b_{1}(X_{m, n})=0$, $b^+(X_{m, n})=4m+2n-1 \equiv 3 \ (\bmod \ 4)$, $b_{1}(Y_{\ell})=0$ and $b^+(Y_{\ell})=3$ holds. Therefore, Theorem \ref{speical-ein} tells us that, $M(m, n, \ell, g, h, \ell_{1})$ cannot admit any Einstein metric if
\begin{eqnarray*}
4(2+\ell_{1} + 1)  \geq \frac{1}{3}\Big( 2\chi(X_{m,n})+3\tau(X_{m,n})+ 2\chi(Y_{\ell})+3\tau(Y_{\ell})+4(1-h)(1-g) \Big),
\end{eqnarray*}
equivalently,
\begin{eqnarray*}
\ell_{1} \geq \frac{1}{12}\Big( 8n+4(1-h)(1-g) \Big) - 3, 
\end{eqnarray*}
where we used $ 2\chi(X_{m,n})+3\tau(X_{m,n})=8n$ and $2\chi(Y_{\ell})+3\tau(Y_{\ell})=0$. However, this inequality is nothing but the inequality (\ref{ein-in-3}) above. Hence, $M(m, n, \ell, g, h, \ell_{1})$ cannot admit any Einstein metric. So each of the topological spin manifolds in (\ref{homeo-spin}) admits infinitely many distinct smooth structures for which no compatible Einstein metric exists. \par
Next we shall prove that each of the topological spin manifolds in (\ref{homeo-spin}), which will be denoted by $M$ for simplicity, has positive volume entropy and satisfies the strict Gromov--Hitchin--Thorpe inequality (\ref{ght}). \par
First of all, by Corollary \ref{main-cor}, we are able to obtain the following bound on the volume entropy of $M$: 
\begin{eqnarray}\label{tel-1}
0< \frac{8}{27{\pi}^2} (g-1)(h-1) \leq \frac{1}{54 {\pi}^2}\lambda(M)^4 \leq \frac{128}{27}(g-1)(h-1).
\end{eqnarray}
In particular, $\lambda(M)>0$ holds. On the other hand, by a direct computation, we are able to get  
\begin{eqnarray}\label{euler-M-1}
2\chi(M)+3\tau(M) = 8n+4(g-1)(h-1)-4(3+\ell_{1}), 
\end{eqnarray}
\begin{eqnarray}\label{euler-M-2}
2\chi(M)-3\tau(M) = 8(12m+n)+96+4(g-1)(h-1)-4(3+\ell_{1}),
\end{eqnarray}
where notice that, for a K3 surface, we have $2\chi+3\tau=0$ and $2\chi-3\tau=96$. \par
Now, by multiplying both sides of (\ref{ein-in-1}) by $4$, we have
\begin{eqnarray*}
8n-\frac{20}{27}(g-1)(h-1) > 4(\ell_{1}+3).
\end{eqnarray*}
Equivalently, 
\begin{eqnarray*}
8n+4(g-1)(h-1) -4(\ell_{1}+3) > \frac{128}{27}(g-1)(h-1).
\end{eqnarray*}
This inequality, (\ref{tel-1}) and (\ref{euler-M-1}) imply
\begin{eqnarray*}
2\chi(M) + 3\tau(M) > \frac{1}{54 {\pi}^2}\lambda(M)^4.
\end{eqnarray*}
Similarly, by multiplying both sides of (\ref{ein-in-2}) by $4$, we get
\begin{eqnarray*}
8(12m+n)-\frac{20}{27} (g-1)(h-1)+84 > 4\ell_{1}.
\end{eqnarray*}
Namely,
\begin{eqnarray*}
8(12m+n)+96+4(g-1)(h-1)-4(3+\ell_{1}) > \frac{128}{27}(g-1)(h-1).
\end{eqnarray*}
This inequality, (\ref{tel-1}) and (\ref{euler-M-2}) tells us that the following holds:
\begin{eqnarray*}
2\chi(M) - 3\tau(M) > \frac{1}{54 {\pi}^2}\lambda(M)^4.
\end{eqnarray*}
Therefore, the spin 4-manifold $M$ satisfies the strict Gromov--Hitchin--Thorpe inequality (\ref{ght}) as desired:
\begin{eqnarray*}
2\chi(M) - 3|\tau(M)| > \frac{1}{54 {\pi}^2}\lambda(M)^4.
\end{eqnarray*}
Hence, the spin 4--manifold $M$ has the desired properties. Because we have infinitely many choices for the above integers $(g, h, m, n, \ell_{1})$, we are able to conclude that there exists an infinite family of topological spin 4--manifolds with the desired properties, as promised.
\end{proof}

Similarly, we are able to prove:

\begin{thm}\label{volume-entropy-2}
There exists an infinite family of topological non-spin four--manifolds satisfying the three properties in Theorem \ref{volume-entropy-1}.
\end{thm}

\begin{proof}
The proof is quite similar to the spin case. For completeness and in order to help the reader, let us include the proof. For any fixed pair $(g, h)$ of odd integers $\geq 3$,  by taking a suitable pair $(m, n)$ of integers satisfying $4m+2n-1 \equiv 3 \ (\bmod \ 4)$, $m \geq 2$ and $n \geq 1$, we are able to find at east one positive integer $\ell_{2}$ satisfying the following three inequalities
\begin{eqnarray}\label{ein-in-12}
8n-\frac{20}{27} (g-1)(h-1)-12 > \ell_{2}.
\end{eqnarray}
\begin{eqnarray}\label{ein-in-22}
8(n+12m) -\frac{20}{27}(g-1)(h-1)+84 > -5\ell_{2}.
\end{eqnarray}
\begin{eqnarray}\label{ein-in-32}
\ell_{2} \geq \frac{1}{3}\Big(8n+4(g-1)(h-1)\Big)-12.
\end{eqnarray}
Notice that we have infinitely many choices of such pairs $(m, n)$ and of $\ell_{2}$. \par
 For any such five integers $(m, n, g, h, \ell_{2})$ and for any new integer $\ell \geq 0$, consider the following connected sum:
\begin{eqnarray*}
N(m, n, \ell, g, h, \ell_{2}):=X_{m, n} \# {Y}_{\ell} \# (\Sigma_{g} \times \Sigma_{h}) \# \ell_{2} \overline{{\mathbb C}{P}^{2}},
\end{eqnarray*}
where $X_{m, n}$ is again Gompf's simply connected symplectic spin 4--manifold and  ${Y}_{\ell}$ is again obtained from $Y_{0}$ by performing a logarithmic transformation of order $2 \ell + 1$ on a non-singular fiber of $Y_{0}$. Notice that $N(m, n, \ell, g, h, \ell_{2})$ is homeomorphic to the following non-spin 4--manifold:
\begin{eqnarray}\label{homeo-non-spin}
X_{m, n} \# K3 \# (\Sigma_{g} \times \Sigma_{h}) \# \ell_{2}\overline{{\mathbb C}{P}^{2}}.
\end{eqnarray}
For any fixed $(m, n, g, h, \ell_{1})$, as in the spin case, we are able to see that the sequence $\{N(m, n, \ell, g, h, \ell_{2}) \ | \ \ell \in {\mathbb N} \}$ contains infinitely many distinct diffeotypes. \par
Moreover, we can also see that $N(m, n, \ell, g, h, \ell_{2})$ cannot admit any Einstein metric. In fact, Theorem \ref{speical-ein} tells us that, for any fixed $(m, n, \ell, g, h, \ell_{2})$, each 4--manifold $N(m, n, \ell, g, h, \ell_{2})$ cannot admit any Einstein metric if
\[
4(2+ 0+ 1) + {\ell}_{2}  \geq \frac{1}{3}\Big( 2\chi(X_{m,n})+3\tau(X_{m,n})+ 2\chi(Y_{\ell})+3\tau(Y_{\ell})+4(1-h)(1-g) \Big), 
\]
namely, 
\begin{eqnarray*}
\ell_{1} \geq \frac{1}{12}\Big( 8n+4(1-h)(1-g) \Big) -12, 
\end{eqnarray*}
where notice again that we have  $b^+(X_{m, n})=4m+2n-1 \equiv 3 \ (\bmod \ 4)$, $b^+(Y_{\ell})=3$, $2\chi(X_{m,n})+3\tau(X_{m,n})=8n$ and $2\chi(Y_{\ell})+3\tau(Y_{\ell})=0$. This inequality is nothing but the inequality (\ref{ein-in-32}) above. Hence, $N(m, n, \ell, g, h, \ell_{2})$ cannot admit any Einstein metric as desired. Hence, the topological non-spin manifold (\ref{homeo-non-spin}) admits infinitely many distinct smooth structures for which no compatible Einstein metric exists. \par
Finally, we shall prove that the manifolds (\ref{homeo-non-spin}), which will be denoted by $N$ for simplicity, have positive volume entropy and satisfy the strict Gromov--Hitchin--Thorpe inequality (\ref{ght}). As in the spin case, Corollary \ref{main-cor} tells us again that the following holds.
\begin{eqnarray}\label{tel-12}
0< \frac{8}{27{\pi}^2} (g-1)(h-1) \leq \frac{1}{54 {\pi}^2}\lambda(N)^4 \leq \frac{128}{27}(g-1)(h-1).
\end{eqnarray}
Hence, we have $\lambda(N)>0$. \par
On the other hand, we get the following by a direct computation:
\begin{eqnarray}\label{euler-M-12}
2\chi(N)+3\tau(N) = 8n+4(g-1)(h-1)-12 -{\ell}_{2}, 
\end{eqnarray}
\begin{eqnarray}\label{euler-M-22}
2\chi(M)-3\tau(M) = 8(12m+n)+84+4(g-1)(h-1)+5\ell_{2}. 
\end{eqnarray}
Now, notice that inequality (\ref{ein-in-12}) is equivalent to 
\begin{eqnarray*}
8n+4(g-1)(h-1)-12-{\ell}_{2} > \frac{128}{27}(g-1)(h-1).
\end{eqnarray*}
This inequality, (\ref{tel-12}) and (\ref{euler-M-12}) imply
\begin{eqnarray*}
2\chi(N) + 3\tau(N) > \frac{1}{54 {\pi}^2}\lambda(N)^4.
\end{eqnarray*}
Similarly, inequality (\ref{ein-in-22}) can be rewritten as
\begin{eqnarray*}
8(12m+n)+4(g-1)(h-1)+84 + 5{\ell}_{2} > \frac{128}{27}(g-1)(h-1).
\end{eqnarray*}
This inequality, (\ref{tel-12}) and (\ref{euler-M-22}) implies
\begin{eqnarray*}
2\chi(N) - 3\tau(N) > \frac{1}{54 {\pi}^2}\lambda(N)^4. 
\end{eqnarray*}
Thus, the non-spin 4--manifold $N$ satisfies the strict Gromov--Hitchin--Thorpe inequality (\ref{ght})
\begin{eqnarray*}
2\chi(N) - 3|\tau(N)| > \frac{1}{54 {\pi}^2}\lambda(N)^4.
\end{eqnarray*}
Therefore, the non--spin 4-manifold $N$ has the desired properties. Since we have infinitely many choices for the above integers $(g, h, m, n, \ell_{2})$, we are able to conclude that there exists an infinite family of topological non-spin 4--manifolds with the desired properties.
\end{proof}

It is clear that Theorem \ref{main-B} now follows from Theorem \ref{volume-entropy-1} or Theorem \ref{volume-entropy-2}.


\providecommand{\bysame}{\leavevmode\hbox to3em{\hrulefill}\thinspace}
\providecommand{\MR}{\relax\ifhmode\unskip\space\fi MR }
\providecommand{\MRhref}[2]{%
  \href{http://www.ams.org/mathscinet-getitem?mr=#1}{#2}
}
\providecommand{\href}[2]{#2}

\end{document}